\documentclass[12pt]{article}
\usepackage{amssymb}
\usepackage{amsthm}
\usepackage{pgf,tikz}
\usepackage{graphicx}
\usepackage{amssymb}
\usepackage{graphics}
\usepackage{amsthm}

\input pictex.tex

\theoremstyle{definition}

\begin{document}

\date{}
\author{Andr\'es Navas\footnote{Funded by Fondecyt Project 1200114}} 

\title{(Un)distorted diffeomorphisms in different regularities}
\maketitle

\vspace{-0.4cm}

\noindent {\bf Abstract.} We build the first examples of diffeomorphisms that are distorted in a group 
of $C^1$ diffeomorphisms yet undistorted in the corresponding group of $C^2$ diffeomorphisms. This 
explicit construction is performed for the closed interval. 


\vspace{0.7cm}

Recall that an element $f$ in a group $G$ is said to be distorted if there exists a finitely generated subgroup 
$\Gamma \subset G$ that contains $f$ such that
$$\lim_{n \to \infty} \frac{\| f^n \|_{_{\Gamma}}}{n} = 0.$$
Here, $\| \cdot \|_{_{\Gamma}}$ stands for the word-length of group elements with respect to any fixed finite generating system of $\Gamma$ 
(the definition of distortion elements above is independent of this choice, since any two such word-lengths are bi-Lipschitz equivalent). It is 
worth pointing out that the limit involved in this definition always exists, due to the (easy to check) subadditive inequality
$$\| f^{m+n} \| \leq \| f^n \| + \| f^m \| $$
and the well-known result commonly called Fekete's Lemma \cite{fekete}.

This notion was introduced by Gromov in \cite{gromov}. Over the last years, it has attracted the attention of people 
not only from group theory and geometry, but also dynamics (see for instance \cite{FH, H}). In this direction, it has been 
successfully used to exhibit obstructions to group actions on  manifolds. Nevertheless, in most of the literature, the fact that 
groups of diffeomorphisms behave differently in distinct regularities is not considered. In this note we would like 
to address the following natural question:

\vspace{0.4cm}

\noindent{\bf Question.} Given a compact manifold $M$ and integers $1 \leq r < s$, does there exist a $C^s$ diffeomorphism 
that is distorted in $\mathrm{Diff}_+^{r} (M)$ yet undistorted in $\mathrm{Diff}_+^s (M)$~?

\vspace{0.4cm}

It is very likely that the answer to this question is affirmative. However, it seems hard to prove it using standard dynamical 
methods, such as genericity type arguments in a well-chosen subset of diffeomorphisms. Namely, as it was pointed out in \cite{CF}, 
diffeomorphisms of positive topological entropy are undistorted in the group of $C^1$ diffeomorphisms (see the Appendix 2 for a 
discussion of this). Therefore, the question above involves zero-entropy diffeomorphisms, which are difficult to tackle. (Relevant 
recent progress on this topic that may be useful in this setting includes \cite{CP,FHandel}.) 

Prototypes of zero-entropy maps are homeomorphisms of 1-manifolds, as it is widely known (see again Appendix 2 on this). 
In this work, we restrict to this framework and we provide a first result in smooth categories concerning the Question above.

\vspace{0.4cm}

\noindent{\bf Main Theorem.} There exist $C^{\infty}$ diffeomorphisms of the closed interval that are $C^1$-distorted yet $C^2$-undistorted.

\vspace{0.4cm}

It is worth stressing that an analogous result holds in lower regularity.     
Indeed, every homeomorphism of the sphere (in any dimension) is $C^0$-distorted, yet it cannot be $C^1$-distorted 
if it is a $C^1$ diffeomorphism that admits hyperbolic periodic points (see \cite{CF}). It is easy to prove a similar  
statement for homeomorphisms of the interval (just use the fact that every nontrivial interval homeomorphism 
is conjugate to its square \cite{Ghys}).  

To prove our result we put together two different ideas. On the one hand, in \S 1, we use the recently introduced notion of {\em asymptotic distortion} 
to ensure that diffeomorphisms in a suitable family are $C^2$ undistorted. On the other hand, in \S 2, we use methods arising in the study of centralizers 
to build finitely many $C^1$ diffeomorphisms that generate groups in which the above examples arise as distorted elements.

The construction that we develop gives us almost sharp examples in relation to the {\em degree of distortion}. On the one hand, 
our $C^{\infty}$ diffeomorphisms $f$ are almost logarithmically distorted, in the sense that $\| f^n \|$ growths slightly faster than logarithmically 
in $n$. On the other hand, no nontrivial $C^1$ diffeomorphism can be more than logarithmically distorted. See \S 3 for the details on this.

It would be certainly worthwhile to improve our methods so that the resulting diffeomorphisms are also distorted in the group of $C^{1+\alpha}$ 
diffeomorphisms for every $\alpha < 1$. (Compare \cite{BKK, DKN}.) This would obviously require more accurate estimates that we plan to 
develop in a complementary (more technical) work.


\section{On the asymptotic distortion of $C^2$ distorted diffeomorphisms}

All maps in this note will be assumed to be orientation preserving.

Given a $C^{1+bv}$ diffeomorphism $f$ of a compact 1-manifold (either the circle or the closed interval), we denote by 
$\mathrm{var} (\log Df)$ the total variation of the logarithm of its derivative. (In this notation, ``$bv$'' stands for ``bounded 
variation''.) The {\em asymptotic distortion} of $f$, denoted $\mathrm{dist}_{\infty} (f)$, is defined as the limit 
$$\lim_{n \to \infty} \frac{\mathrm{var} (\log Df^n )}{n}.$$
This notion was introduced in \cite{duke} with a focus on circle diffeomorphisms, and later studied in \cite{EN} for diffeomorphisms of the 
interval. The relation with distorted elements is clear: if $f$ is distorted in the group of $C^{1+bv}$ diffeomorphisms, then its asymptotic 
distortion must vanish (this directly follows from the subaditivity property of $\mathrm{var} (\log D (\cdot))$; see \cite[Corollary 4]{EN}).  

A nice family of diffeomorphisms with positive asymptotic distortion may be described as follows: Start with a $C^{1+bv}$ diffeomorphism of 
$[0,1]$ with vanishing asymptotic distortion and no fixed point in $]0,1[$. Let $I$ be a {\em fundamental domain} for the action of $f$, that is, 
an open interval with endpoints $x_0, f (x_0)$ for a certain $x_0 \in \, ]0,1[$. Let $g$ be any nontrivial $C^{1+bv}$ diffeomorphism of $]0,1[$ 
supported on $I$. Then the diffeomorphism $\bar{f} := fg$ is undistorted in $\mathrm{Diff}_+^{1+bv} ([0,1])$, hence in $\mathrm{Diff}_+^2([0,1])$; 
see \cite[Lemmas 2.2 and 7.2]{EN}.

\begin{center}
\begin{tikzpicture}[x=1.0cm,y=1.0cm]
\draw [line width=1.pt] (0.,2.3)-- (0.,1.7);
\draw [line width=1.pt] (0.,2.)-- (9.,2.);
\draw [line width=1.pt] (9.,2.3)-- (9.,1.7);
\draw [line width=1.pt] (2.4,2.2)-- (2.4,1.8);
\draw [line width=1.pt] (3.8,2.2)-- (3.8,1.8);
\draw [line width=1.pt] (5.2,2.2)-- (5.2,1.8);
\draw [line width=1.pt] (6.6,2.2)-- (6.6,1.8);
\draw [line width=1.pt] (0.,1.0)-- (0.,0.6);
\draw [line width=1.pt] (3.8,1.0)-- (3.8,0.6);
\draw [line width=1.pt] (5.2,1.0)-- (5.2,0.6);
\draw [line width=1.pt] (9.,1.0)-- (9.,0.6);
\draw [line width=1.pt] (5.2,0.8)-- (6.3,0.8);
\draw [line width=1.pt] (7.9,0.8)-- (9.,0.8);
\draw [line width=1.pt] (2.7,0.8)-- (3.8,0.8);
\draw [line width=1.pt] (0.,0.8)-- (1.1,0.8);
\draw (2.0,2.8) node[anchor=north west] {$x_{-1}$};
\draw (3.5,2.8) node[anchor=north west] {$x_0$};
\draw (4.9,2.8) node[anchor=north west] {$x_1$};
\draw (6.3,2.8) node[anchor=north west] {$x_2$};
\draw (7.41,2.65) node[anchor=north west] {\scalebox{1.2}{$.\,.\,.$}};
\draw (0.656,2.65) node[anchor=north west] {\scalebox{1.2}{$.\,.\,.$}};
\draw (1.15,1.15) node[anchor=north west] {$g=Id$};
\draw (6.35,1.15) node[anchor=north west] {$g=Id$};
\draw (4.3,0.5) node[anchor=north west] {$g$};
\draw (-0.25,1.75) node[anchor=north west] {$0$};
\draw (8.75,1.75) node[anchor=north west] {$1$};
\draw (5.0,3.9) node[anchor=north west] {$f$};
\draw [shift={(3.7,2.64)},line width=1.pt]  plot[domain=-0.07:3.22,variable=\t]({1.*0.60*cos(\t r)+0.*0.60*sin(\t r)},{0.*0.60*cos(\t r)+1.*0.60*sin(\t r)});
\draw [shift={(5.2,2.64)},line width=1.pt]  plot[domain=-0.07:3.22,variable=\t]({1.*0.60*cos(\t r)+0.*0.60*sin(\t r)},{0.*0.60*cos(\t r)+1.*0.60*sin(\t r)});
\draw [shift={(4.5,0.9)},line width=1.pt]  plot[domain=-3.14:0.,variable=\t]({1.*0.35*cos(\t r)+0.*0.35*sin(\t r)},{0.*0.35*cos(\t r)+1.*0.35*sin(\t r)});
\draw [line width=1.pt] (4.85,0.9)-- (4.85,1.2);
\draw [line width=1.pt] (4.15,0.9)-- (4.15,1.1);
\draw (4.61,1.5) node[anchor=north west] {\scalebox{.8}{$\blacktriangle$}};
\draw (4.05,2.75) node[anchor=north west] {\scalebox{.8}{$\blacktriangledown$}};
\draw (5.55,2.75) node[anchor=north west] {\scalebox{.8}{$\blacktriangledown$}};
\end{tikzpicture}
\end{center}

If $f$ is $C^{\infty}$ tangent to the identity at the endpoints, we may extend it to a larger interval (say, $[-1,2]$) by the identity 
outside $[0,1]$. It follows from the definition that this new diffeomorphism still has zero asymptotic distortion. Similarly, the 
extension of $\bar{f}$ has nonzero asymptotic distortion, and it is hence undistorted in $\mathrm{Diff}^2_+([-1,2])$.

\vspace{0.35cm}

\noindent{\bf Remark.} Starting with the seminal work of Rosendal \cite{ros}, over the last years there has been a lot of interest in understanding 
the geometry of (big) topological groups. In this setting, the word-metric is considered with respect to neighborhoods of the identity. The notion 
of distorted element that arises this way is weaker that the one considered in this paper, hence that of undistorted element is stronger. However, 
the proof above shows that $\bar{f}$ is still undistorted in this more restricted sense.


\section{How to force $C^1$ distortion}

We next explain how to perform the construction above so that the resulting diffeomorphism $\bar{f}$ of $[-1,2]$ is 
$C^1$ distorted. This will be accomplished in several steps.

\vspace{0.5cm}

\noindent{\bf {\em Building the starting diffeomorphism $f$.}} 
Start with the vector fields $\hat{\mathcal{X}}$ and $\mathcal{X}$ on the real line 
whose time-1 maps are, respectively, 
$$\hat{F} 
:= \hat{\mathcal{X}}^1
: x \mapsto 2x \qquad { \mathrm{and} } \qquad F
:= \mathcal{X}^1
: x \mapsto x+1.$$
It is well known that there exists a $C^{\infty}$ diffeomorphism $\varphi: \mathbb{R} \to ]0,1[$ such that  \, 
$\hat{\mathcal{Y}} := \varphi_* (\hat{\mathcal{X}})$ \, and \, $\mathcal{Y} := \varphi_* (\mathcal{X})$ \, 
extend to the endpoints of $[0,1]$ as infinitely flat vector fields. (See for instance \cite{navas-book, tsu}.) 
Let us denote $\hat{f} := \hat{\mathcal{Y}}^1$ and $f := \mathcal{Y}^1$, 
which we view as diffeomorphisms of $[-1,2]$ that coincide with the identity outside $[0,1]$. It follows from \cite{EN} that $f$ has 
vanishing asymptotic distortion. Actually, $f$ is distorted in the group $\langle \hat{f}, f \rangle$. Indeed, the (Baumslag-Solitar) 
relation $\hat{f} f \hat{f}^{-1} = f^2$ implies that, for each $k \geq 1$, 
$$\hat{f}^i f \hat{f}^{-i} = f^{2^i}.$$
As a consequence, $\| f^n \|$ is of order at most $O ( \log(n))$. (Actually, it is of order $O (\log(n))$; this is very well known because 
of the Baumslag-Solitar group structure, but also follows from \S 3 below.)

Let $x_0 := \varphi (0)$ and, for each $n \in \mathbb{Z}$, let $x_n := f^n (x_0)$. 
Notice that $x_n = \varphi (n)$. Denote also $x_{-1/2} := \varphi (-1/2)$ and $x_{-3/4} := \varphi (-3/4)$.

\vspace{0.5cm}

\noindent{\bf {\em Building the diffeomorphism $g$.}} 
Let us denote by $\varphi_0$ the affine diffeomorphism sending $[x_0,x_1]$ onto $[0,1]$. Let $g := \varphi_0^{-1} f \varphi_0$. 
This can be extended to $[-1,2]$ by the identity outside $[x_0, x_1]$. As above, we let $\bar{f} := fg$.

\vspace{0.5cm}

\noindent{\bf {\em A first estimate to prove distortion.}} 
Notice that 
$$\bar{f}^n f =  (fg)^n f = f^{n+1} (f^{-n} g f^{n}) (f^{-(n-1)} g f^{n-1}) \cdots (f^{-2} g f^2)  (f^{-1} g f),$$
hence
$$f^{-n} (f^{-1} \bar{f}^n f) =  (f^{-n} g f^{n}) (f^{-(n-1)} g f^{n-1}) \cdots (f^{-2} g f^2)  (f^{-1} g f).$$
Let us denote by $h_n$ this diffeomorphism. Clearly, it is supported on $[x_{-n}, x_0]$. 
Moreover, since $f$ is distorted in $\langle \hat {f}, f  \rangle$, 
if we show that $\| h_n \|$ is $o(n)$ along an increasing (infinite) sequence of integers, then the same will 
be true for $\| f^{-1} \bar{f}^n f \|$, hence for $\| \bar{f}^n \|$.  
This will imply that $\bar{f}$ is a distorted element.

\vspace{0.5cm}

\noindent{\bf {\em Building two auxiliary diffeomorphisms.}} 
We define two diffeomorphisms $\hat{h}$ and $h$ as follows:

\vspace{0.18cm}

\noindent (i) They act by the identity outside $[0,1]$.

\vspace{0.18cm}

\noindent (ii) On each interval $f^n ([x_0,x_1])$, the diffeomorphism $\hat{h}$ (resp. $h$) coincides with the 
$s_n$-time map (resp. $t_n$-time map) of the flow of the vector field $f^n_* \big(\varphi_0^* (\hat{\mathcal{Y}} ) \big)$ 
(resp. $f^n_* \big(\varphi_0^* (\mathcal{Y})) $\big). 

\vspace{0.18cm}

Here, $s_n$ and $t_n$ are sequences of real numbers to be defined. 
According to \cite[Lemma 4.1]{DF}\footnote{Related results appear in \cite{BF,Fa}. For the reader's convenience, 
in the Appendix 1, we give a direct and simple proof that does not rely on Kopell's (difficult) estimates.}, 
in order to ensure that $\hat{h}$ and $h$ are $C^1$ at the endpoints 
(and hence $C^1$ diffeomorphisms), it is necessary and sufficient to impose the conditions $s_n \rightarrow 0$ 
and $t_n \rightarrow 0$ as $n \rightarrow \pm \infty$. In practice, we will choose these sequences as follows: 

\vspace{0.18cm} 

\noindent (iii) If $2^{i-1} \leq n < 2^{i}$ for a certain positive even integer $i$, we let 
$$s_n := \log_2 \left( 1 - \frac{1}{\sqrt{\ell_{i/2}}} \right) \qquad \mbox{ and } \qquad t_n := \frac{1}{ \sqrt{\ell_{i/2}} },$$
where $\ell_j$ is any prescribed sequence of positive integers diverging to infinity.

\vspace{0.18cm} 

\noindent (iv) For any other value of $n$, we let $s_n = t_n := 0$.

\vspace{0.5cm}

\noindent{\bf {\em An extra diffeomorphism to control supports.}} 
Finally, we let $\psi$ be a $C^{\infty}$ diffeomorphism of $[-1,2]$ such that:

\vspace{0.18cm}

\noindent (v) $\psi$ coincides with the identity on $[x_{-1/2}, x_0]$.

\vspace{0.18cm}

\noindent (vi) $\psi (x_{-3/4}) = 0$ and $\psi (x_1) = 1$.

\vspace{0.18cm}

\noindent Computations below are made in the group $\Gamma := \langle \hat{f}, f, g, \hat{h}, h, \psi \rangle$ (contained in $\mathrm{Diff}^1_+([-1,2])$). 

\vspace{0.5cm}

\noindent{\bf {\em On the conjugates of $h$ and $\hat{h}$.}} Fix an integer $n$ of the form $2^i$, where $i$ is a large even number. 
We start by considering the conjugate maps $f^{-n} h f^{n}$, later $\hat{f}^{-i} f^{-n} h f^n \hat{f}^{i}$, and finally $\psi \hat{f}^{-i} f^{-n} h f^n \hat{f}^{i} \psi^{-1}$. 
Since 
$$\mathrm{supp} (h) \subset [0, x_{2^{i-2}}] \cup [x_{2^{i-1}}, x_{2^{i}}] \cup [x_{2^{i+1}},1],$$
we have
$$\mathrm{supp} (f^{-n} h f^{n}) \subset [0, x_{2^{i-2} - 2^i}] \cup [x_{2^{i-1} - 2^i}, x_{2^i - 2^i}] \cup [x_{2^{i+1} - 2^i},1],$$
 that is,
 $$\mathrm{supp} (f^{-n} h f^{n}) \subset [0, x_{2^{i-2} - 2^i}] \cup [x_{-2^{i-1}},x_0] \cup [x_{2^i},1].$$
Since \, $(2^{i-2} - 2^i) / 2^i = -3 / 4$ \, and \, $- 2^{i-1} / 2^i = -1/2$, \, this implies
$$\mathrm{supp} (\hat{f}^{-i} f^{-n} h f^n \hat{f}^{i} ) \subset [0,x_{-3/4}] \cup [x_{-1/2},x_0] \cup [x_1,1].$$
Finally, the properties (v) and (vi) imposed to $\psi$ imply  
$$\mathrm{supp} ( \psi \hat{f}^{-i} f^{-n} h f^n \hat{f}^{i} \psi^{-1} ) \subset [-1,0] \cup [x_{-1/2},x_0] \cup [1,2].$$ 

Similarly,
$$\mathrm{supp} (\hat{f}^{-i} f^{-n} \hat{h} f^n \hat{f}^{i} ) \subset [0,x_{-3/4}] \cup [x_{-1/2},x_0] \cup [x_1,1],$$
and
$$\mathrm{supp} ( \psi \hat{f}^{-i} f^{-n} \hat{h} f^n \hat{f}^{i} \psi^{-1} ) \subset [-1,0] \cup [x_{-1/2},x_0] \cup [1,2].$$

\vspace{0.5cm}

\noindent{\bf {\em A commutator of conjugates of $h$ and $\hat{h}$.}} Let us denote \, $a_n := \hat{f}^{-i} f^{-n} h f^n \hat{f}^{i}$ \, and \, 
$b_n~:=\psi \hat{f}^{-i} f^{-n} \hat{h} f^n \hat{f}^{i} \psi^{-1}.$ \, We next compute $c_n := [a_n, b_n] = a_n (b_n a_n b_n^{-1})^{-1}$. 
We notice that:

\vspace{0.18cm}

\noindent -- The intervals $[-1,0]$ and $[1,2]$ are invariant under both $a_n$ and $b_n$; besides, $a_n$ equals the identity therein. 
Therefore, on $[-1,0] \cup [1,2]$, the map $c_n$ acts trivially.

\vspace{0.18cm}

\noindent -- The intervals $[0, x_{-1/2}]$ and $[x_0,1]$ are invariant under both $a_n$ and $b_n$; besides, $b_n$ equals the identity therein. 
Therefore, on $[0, x_{-1/2}] \cup [x_0, 1]$, the map $c_n$ acts trivially as well.

\vspace{0.18cm}

\noindent -- By property (v) above, the restriction of $b_n$ to $[x_{-1/2},x_0]$ coincides with that of $\hat{f}^{-i} f^{-n} \hat{h} f^n \hat{f}^{i}$. 
Therefore,  the restriction of $c_n$ to this interval equals that of 
$$\hat{f}^{-i} f^{-n} (h \hat{h} h^{-1} \hat{h}^{-1}) f^n \hat{f}^{i}.$$

\vspace{0.5cm}

\noindent{\bf {\em Retrieving $h_n$ as a power of a commutator.}} Recall that  
$$h_n := f^{-n} (f^{-1} \bar{f}^n f) =  (f^{-n} g f^{n}) (f^{-(n-1)} g f^{n-1}) \cdots (f^{-2} g f^2)  (f^{-1} g f).$$
We claim that, for $n = 2^i$,  
$$h_{n/2}= \big( \hat{f}^i c_n \hat{f}^{-i} \big)^{\, \ell_{i/2}}.$$
Indeed, by the discussion above, $\hat{f}^i c_n \hat{f}^{-i}$ is supported on $\hat{f}^i ([x_{-1/2}, x_0]) = [x_{-2^{i-1}}, x_0]$, 
which coincides with the support of $h_{n/2} = h_{2^{i-1}}$. Moreover, 
$$\hat{f}^i c_n \hat{f}^{-i} = ( f^{-n} h f^n ) \Big[ (f^{-n} \hat{h} f^n) ( f^{-n} h^{-1} f^n ) ( f^{-n} \hat{h}^{-1} f^n ) \Big].$$
Furthermore, by (ii) and (iii), on each interval $[x_{k}, x_{k+1}]$, with $-2^{i-1} \leq k < 0$, the map 
$\hat{h}$ (resp. $h$) arises as the time-$s$ (resp. time-$t$) map of the flow of $f^k_* (\varphi_0^* ( \hat{\mathcal{Y}}))$ 
(resp. $f^k_* ( \varphi_0^* (\mathcal{Y}))$), where 
$$s = \log_2 \left( 1 - \frac{1}{\sqrt{\ell_{i/2}}} \right) \qquad \mbox{ and } \qquad t = \frac{1}{ \sqrt{\ell_{i/2}} }.$$
The Lie-algebra relation between $\hat{\mathcal{Y}}$ and $\mathcal{Y}$ then implies that, on each such interval $[x_k,x_{k+1}]$, 
the commutator \, $\hat{f}^i c_n \hat{f}^{-i} = ( f^{-n} h f^n ) \big[ (f^{-n} \hat{h} f^n) ( f^{-n} h^{-1} f^n ) ( f^{-n} \hat{h}^{-1} f^n ) \big]$ \, 
coincides with the time-$r$ map of the flow of $f^k_* (\mathcal{Y})$, where 
$$r = (1-2^s) \, t = \frac{1}{\ell_{i/2}}.$$
Since the restriction of $h_n$ to $[x_k,x_{k+1}]$ is nothing but the time-1 map of the flow of $f^k_* (\mathcal{Y})$, this implies our claim.

\vspace{0.5cm}

\noindent{\bf {\em A final estimate to prove distortion.}} Since\, $h_{n/2} = \big( \hat{f}^i c_n \hat{f}^{-i} \big)^{\, \ell_{i / 2}}$, we have 
$$\| h_{n/2} \| \leq \ell_{i/2} \, \| \hat{f}^i c_n \hat{f}^{-i}  \| \leq \ell_{i/2} \, [ 2 i + \| c_n \| ] \leq 2 \, \ell_{i/2} \, [ i + \| a_n \| + \| b_n \| ].$$
Using the definitions of $a_n = \hat{f}^{-i} f^{-n} h f^n \hat{f}^{i}$ and $b_n = \psi \hat{f}^{-i} f^{-n} \hat{h} f^n \hat{f}^{i} \psi^{-1}$, 
this yields 
$$\| h_{n/2} \| \leq 2 \, \ell_{i/2} \, [ i +  (2i + 1 + 2 \|f^n\| ) + (2i + 3 + 2 \|f^n\|) ] = 2 \, \ell_{i/2} \, [5i+4+4\|f^n\|].$$
Since $i \sim \log(n)$ and $\| f^n \| = O (\log(n))$, we conclude that 
$$\| h_{n/2} \| = O \big( \ell_{i/2} \log (n) \big).$$
Therefore, if we choose the sequence $\ell_j$ so that $\ell_{i/2} \log (n) = o(n)$, we have that $\| h_n \| = o (n)$. As previously 
discussed, this implies that $\bar{f}$ is a distorted element of $\Gamma$.

\vspace{0.5cm}

\noindent{\bf Remark.} It seems hard to completely describe the group $\Gamma$ above in algebraic terms. In particular, we do not 
know whether it contains free subgroups in two generators. However, it is easy to find elements for which Lemma 3.1 of \cite{KKL} 
may be applied, thus showing that $\Gamma$ contains (many) copies of Thompson's group $F$; in particular, it is not solvable. Whether 
or not $C^2$-undistorted diffeomorphisms of $[0,1]$ may be distorted inside a solvable group of $C^1$ diffeomorphisms seems to be an 
interesting question. (Compare \cite{BMNR}.)


\section{No super-exponential distortion}

Being distorted means that the word-length of powers grows sublinearly, but one may ask whether this growth is even weaker. In the explicit examples 
above, the computations show that $\bar{f}$ can be built as being slightly weaker than exponentially distorted. More precisely, for any sequence $k_n$ of 
positive integers diverging to infinity, $ \| \bar{f}^n \|$ was showed to be of order $\, o (k_n \, \log (n)) \, $ inside an appropriate finitely generated group 
(Actually, this was proved along a subsequence of integers $n$ of the form $2^i$, but it is easy to ensure such a global behavior from this.) 
This is almost sharp, as proven next.

\vspace{0.4cm}

\noindent{\bf Proposition.} {\em Given a nontrivial $C^1$ diffeomorphism $\bar{f}$ inside a finitely generated group of diffeomorphisms 
of a closed interval, there exists a constant $A > 0$ such that, for infinitely many integers $n \geq 0$,}
$$\| \bar{f}^n \| \geq A \log (n).$$

\vspace{0.2cm}

\noindent{\bf Proof.} This is essentially a restatement of an elementary observation from \cite[pages 192-193]{polt-sodin} that claims the following: If $\bar{f}$ 
is a nontrivial diffeomorphism of a closed interval, then there exists an increasing sequence of integers $n_k$ and points $x_k$ such that 
$$D\bar{f}^{n_k} (x_k) \geq n_k.$$
Now, let $\bar{f}$ lie in a group of $C^1$ diffeomorphisms of the interval generated by finitely many elements $\{f_1, \ldots,f_{\ell}\}$, and 
let $C$ be the maximum possible value of all the $D f_i^{\pm 1}$. If, for a fixed $k \geq 1$, we denote $m = m_k := \| \bar{f}^{n_k} \|$, then 
$$\bar{f}^{n_k} = f_{i_1} \cdots f_{i_m}$$
for certain $f_{i_j} \in \langle f_1, \ldots, f_{\ell} \rangle$. This obviously implies that $D \bar{f}^{n_k} (x) \leq C^m$ for all points $x$. 
Choosing $x = x_k$ we conclude that $n_k \leq C^m$, which implies $m \geq \log (n_k) / \log (C)$, that is, 
$$\| \bar{f}^{n_k} \| \geq \frac{\log (n_k)}{ \log (C)},$$
as claimed. 
$\hfill\square$


\section{A remark on (un)distorted circle diffeomorphisms}

By gluing the endpoints of the interval, the previous example can be realized as a 
$C^{\infty}$ circle diffeomorphism with a single global fixed point. More generally, for any rational number 
$p/q$ (mod. $1$), one can adapt the construction to build examples of $C^{\infty}$ circle diffeomorphisms 
of rotation number $p/q$ that are $C^1$ distorted yet $C^2$ undistorted. The situation for irrational 
rotation number is, however, less clear. 

On the one hand, every $C^2$ (actually, $C^{1+bv}$) circle diffeomorphism $f$ of irrational rotation 
number $\rho$ is $C^1$ distorted. This follows as a combination of a slight modification of a theorem 
of Avila (indeed, his method of proof yields $C^1$ distortion provided $C^1$ recurrence is ensured) and a theorem 
of Herman \cite{herman} according to which these diffeomorphisms are recurrent, in the sense that $f^{q_n}$ converges to the 
identity in the $C^1$ topology as $n$ goes to infinity, where $p_n / q_n$ is the sequence of (best) rational approximations of $\rho$ 
(see \cite{NT} for a much simpler proof of this result).

On the other hand, we ignore whether a $C^{\infty}$ (or even $C^2$) diffeomorphism $f$ as above can be undistorted. In 
one direction, the main idea used 
in the construction of $C^2$ undistorted interval diffeomorphisms fails here, since the asymptotic distortion of such an $f$ always vanishes 
\cite{duke}. In the other direction, Avila's method for showing distortion fails not only because of the absence of recurrence in the $C^{2}$ topology, 
but also because another important tool of his proof, namely Mather's simplicity type results, is unavailable in the $C^{2}$ setting. However, it 
is worth mentioning that  one can build $C^{1+bv}$ circle diffeomorphisms of irrational rotation number for which the 
asymptotic distortion is nonzero; these are hence undistorted in $\mathrm{Diff}_+^{1+bv} (\mathrm{S}^1)$ 
(see \cite[Remark 2]{duke} for the details). 


\section{Appendix 1: A generalized version of Kopell's lemma} 
	
\vspace{0.4cm}

In \S 2, we strongly used the result below. For $f$ of class $C^2$, this is Lemma 4.1 in \cite{DF}. Stated for $f$ of class $C^{1+Lip}$, 
this is Exercise 4.1.5 in \cite{navas-book}. Let us stress that we do not know whether this still holds for $f$ merely of class $C^{1+bv}$, 
which is the most general setting for the classical Kopell's lemma.

\vspace{0.4cm}

\noindent{\bf Lemma.} {\em Let $f$ be a $C^{1+Lip}$ diffeomorphism of $[0,1)$ such that $f (x) < x$ for all $x \in (0,1)$. 
Let $x_0$ be a point in $(0,1)$, and denote $x_n := f^n (x_0)$. For each sequence $(g_n)$ of $C^1$ diffeomorphisms 
of $[x_1,x_0]$ that are tangent to the identity at the endpoints, the following conditions are equivalent:}

\vspace{0.1cm}

\noindent (i) {\em The sequence $(g_n)$ converges to the identity in the $C^1$ topology.}

\vspace{0.1cm}

\noindent (ii) {\em The diffeomorphism $g \!: (0,x_0] \to (0,x_0]$ whose restriction to $[x_{n+1},x_n]$ coincides with 
$f^n g_n f^{-n}$ extends to a $C^1$ diffeomorphism of $\, [0,x_0]$ by letting $g(0) = 0$.}

\vspace{0.4cm}

\noindent{\bf Proof.} Assuming (i), let us prove (ii). We need to show that, as $n$ goes to infinity, the derivative $Dg (y)$ 
uniformly converges to $1$ for $y \in [x_{n+1}, x_n]$. Using the relation $g = f^n g_n f^{-n}$ on $[x_{n+1}, x_n]$, we compute 
\begin{equation}\label{cero0}
D g (y) 
= \frac{Df^n (g_n f^{-n} (y))}{Df^n (f^{-n} (y)) } \cdot Dg_n (f^{-n} (y)) 
= \frac{Df^n (z_n)}{Df^n (y_n)} \cdot Dg_n (y_n),
\end{equation}
where $y_n := f^{-n} (y)$ and $z_n := g_n (y_n)$ both belong to $[x_1, x_0]$. Since $g_n$ converges to $Id$ in the $C^1$ 
topology, the factor $Dg_n (y_n)$ converges to $1$ as $n$ goes to infinity. Concerning the quotient of derivatives of $Df^n$, 
we first notice that, given $\varepsilon > 0$, for $n$ large enough we have 
\begin{equation}\label{prim1}
\frac{|z_n - y_n|}{|x_0 - x_1|} 
\leq \frac{| z_n - y_n|}{|y_n - x_1|} 
= \left| \frac{g_n (y_n) - g_n (x_1)}{y_n - x_1} - 1 \right| 
= | Dg_n (\xi_n) - 1| \leq \varepsilon,
\end{equation}
where, in the previous computation, $\xi_n$ is a certain point in $[x_1,x_0]$. Now, for all $i \geq 0$, we have 
\begin{equation}\label{seg2}
\frac{|f^i (z_n) - f^i (y_n)|}{|x_n - x_{n+1}|} = \frac{Df^i (q) \cdot (z_n - y_n)}{Df^i (p) \cdot (x_0 - x_1)}
\end{equation}
for certain points $q = q_{i,n}$ and $p = p_{i,n}$ in $[x_1,x_0]$. Moreover, letting 
$C$ be the Lipschitz constant of $\log (Df)$ on $[0,x_0]$, we have, for all $i \geq 0$,
\begin{eqnarray*}
\left| \log \left( \frac{Df^i (q)}{Df^i (p)} \right) \right| 
&=& \left| \sum_{j=0}^{i-1} \big[ \log Df (f^j (q)) - \log Df (f^j (p)) \big] \right| \\
&\leq& \sum_{j=0}^{i-1} \big| \log Df (f^j (q)) - \log Df (f^j (p)) \big| \\
&\leq& \sum_{j=0}^{i-1} C \, \big| f^j (q) - f^j (p) \big| \\
&\leq& C.
\end{eqnarray*}
According to (\ref{prim1}) and (\ref{seg2}), for a large-enough $n$, this implies 
$$\frac{|f^i (z_n) - f^i (y_n)|}{|x_n - x_{n+1}|} \leq \varepsilon \, e^C.$$
Therefore, 
\begin{eqnarray*} 
\left| \log \left( \frac{Df^n (z_n)}{Df^n (y_n)} \right) \right| 
&=& \left| \sum_{i=0}^{n-1} \big[ \log Df (f^i (z_n)) - \log Df (f^i (y_n)) \big] \right| \\
&\leq& \sum_{i=0}^{n-1} \big| \log Df (f^i (z_n)) - \log Df (f^i (y_n)) \big| \\
&\leq& \sum_{i=0}^{n-1} C \, \big| f^i(z_n) - f^i(y_n) \big| \\
&\leq& C \, \varepsilon \, e^C \sum_{i=0}^{n-1} |x_n - x_{n+1}| \\
&\leq& C \, \varepsilon \, e^C.
\end{eqnarray*}
This implies that the value of the quotient in the expression (\ref{cero0}) lies in $[e^{-\varepsilon\, C'}, e^{\varepsilon\, C'}]$ 
for $C' := Ce^C$. Since $\varepsilon > 0$ was arbitrary, this shows that, as $n$ goes to infinity, the value of $Dg(y)$ 
uniformly converges to $1$ for $y \in [x_{n+1}, x_n]$.

The proof that (ii) implies (i) roughly consists in reversing the arguments above. To begin with, just notice that 
$Dg (0)$ must be equal to $1$, since $g$ fixes all the points $x_n$, which converge to the origin. 
The rest of the computations are straightforward and may be left to the reader. 
$\hfill\square$


\section{Appendix 2: A remark on entropy and distortion} 
	
\vspace{0.4cm}

As it is cleverly noticed in \cite{CF}, an application of Ruelle's inequality proves that $C^1$ diffeomorphisms of compact manifolds with 
positive topological entropy are undistorted. Here we extend this fact to bi-Lipschitz homeomorphisms. 

\vspace{0.45cm}

\noindent{\bf Proposition.} {\em Every bi-Lipschitz homeomorphism of a compact manifold with 
positive topological entropy is undistorted (in the corresponding group of bi-Lipschitz homeomorphisms).}

\vspace{0.45cm}

\noindent{\bf Proof.} This relies on the Bowen-Kouchnirenko inequality \cite{bowen,kouch} for a $C^1$ 
diffeomorphism $f$ of a $d$-dimensional compact manifold:
$$h_{\mathrm{top}} (f) \leq d \cdot \max_{x} \log \big( \| Df (x) ) \| \big).$$
Indeed, this inequality holds (with the very same proof) for bi-Lipschitz homeomorphisms. In this setting, the ``$\max$'' 
should be replaced by an ``$\mathrm{essup}$'', which is controlled by the Lipschitz constant of the map\footnote{This more general version allows reproving 
the well-known fact that homeomorphisms of (compact) 1-manifolds have zero entropy. Indeed, an elementary argument given in \cite{rap} 
establishes that these homeomorphisms are topologically conjugate to bi-Lipschitz homeomorphisms with Lipschitz constant as close to $1$ as desired.}:
$$h_{\mathrm{top}} (f) \leq d \cdot \log ( Lip(f) ).$$ 

Now, let $f$ be a bi-Lipschitz homeomorphism of a compact $d$-dimensional manifold that is distorted, and let 
$\Gamma$ 
be a finite generating system of a group of bi-Lipschitz homeomorphisms containing $f$ such that
$$\lim_{n \to \infty} \frac{\| f^n \|_{\Gamma}}{n} = 0.$$
If $\ell_n := \| f^n \|_{\Gamma}$, then there exist $f_{i_1}, \ldots, f_{i_{\ell_n}}$ in $\Gamma$ 
such that $f^n = f_{i_1} \cdots f_{i_{\ell_n}}$. If we denote by $L$ the maximum of the Lipschitz constants of the maps in $\Gamma$, 
this implies \, $Lip (f^n) \leq L^{\ell_n}.$ \, Thus, 
$$\log (Lip (f^n)) \leq \ell_n \log(L).$$ 
By the Bowen-Kournichenko inequality, this yields
$$h_{\mathrm{top}} (f) = \frac{h_{\mathrm{top}}(f^n)}{n} \leq \frac{1}{n} \, d \cdot \ell_n \log (L),$$
which implies that $h_{\mathrm{top}} (f) = 0$ by letting $n$ go to infinity.
$\hfill\square$

\vspace{0.35cm}

\noindent{\bf Acknowledgments.} I wish to thank Adrien Le Boudec, Sang-Hyun Kim and Christian Rosendal for their 
useful remarks, Jairo Bochi and Godofredo Iommi for a couple of references, and the referee for her/his many corrections.
	

\begin{small}


\vspace{0.3cm}

\noindent Andr\'es Navas

\noindent Dpto. de Matem\'atica y C.C., Univ. de Santiago de Chile (USACH)

\noindent Alameda 3363, Estaci\'on Central, Santiago, Chile

\noindent Email: andres.navas@usach.cl

\end{small}

\end{document}